\documentclass{amsart}

\usepackage{amsmath}
\usepackage{amsfonts}
\usepackage{amsthm}

\newtheorem{definition}{Definition}
\newtheorem{proposition}{Proposition}
\newtheorem{theorem}{Theorem}
\newtheorem{lemma}{Lemma}
\newtheorem{remark}{Remark}

\newcommand{\rep}{\mathop{Rep}}

\begin{document}



\title{Unbounded representations of $q$-deformation of Cuntz algebra}

\author{Vasyl {Ostrovskyi}%
}\email{vo@imath.kiev.ua}
\address{Institute of Mathematics,  National Academy of
Sciences of Ukraine}

\author{Daniil {Proskurin}%
}\email{prosk@univ.kiev.ua}
\address{Kyiv National Taras Shevchenko University, Faculty of
Cybernetics}

\author{Lyudmila {Turowska}%
}\email{turowska@math.chalmers.se}
\address{Chalmers University of Technology, Department of
Mathematics}

\dedicatory{To the memory of Leonid Vaksman}



\begin{abstract}
We study a deformation of the Cuntz-Toeplitz $C^*$-algebra
determined by the relations $a_i^*a_i=1+q a_ia_i^*,\ a_i^*a_j=0$.
We define  well-behaved unbounded $*$-representations of the
$*$-algebra defined by relations above and
classify all such irreducible representations up to
unitary equivalence.
\end{abstract}




\maketitle

\section*{Introduction}
Many of the structures that have been studied recently arise as
deformations of classical objects, e.g. deformations of the
canonical commutation relations (CCR) and the canonical
anti-commutation relations (CAR), quantum groups, quantum
homogeneous spaces, non-commutative  probability  etc (see \cite{
biedenharn, bsp1, gre91, mac89, pus_wor89, vaksman, vaksman1}). From
the physical point of view  important classes of objects come from
the Fock space formalism forming algebras generated by raising and
lowering operators and their numerous generalizations. Examples of
such generalizations include $q$-deformed quantum oscillator algebra
(\cite{biedenharn, gre91, mac89}), twisted CCR (\cite{pus_wor89}),
generalised deformed oscillator (\cite{daskaloyanis}) and more
general quadratic algebras with Wick ordering
(\cite{jor_sch_wer95}). Here unbounded representations arise
naturally since most of the physical observables can not be realized
by bounded operators. During the past 30 years many works concerning
(topological) algebras of unbounded operators and their physical
applications appeared in the literature (see e.g.
\cite{bagarello,inoue,jurzak,sch90} and references therein).

One of the most known $q$-deformations of CCR is $q$-CCR with one
degree of freedom, see \cite{biedenharn, mac89}, i.e. the
$*$-algebra which is generated by elements $a$, $a^*$ satisfying the
commutation relation
\begin{equation}
\label{q1} a^* a=1+q a a^*,
\end{equation}
where  $q\in [0,1)$.

For many degree of freedom there exist several versions of $q$-CCR
algebras, see  \cite{bsp,gre91,marcinek, pus_wor89}. In this paper
we consider a subclass of $q_{ij}$-CCR algebras introduced in
\cite{bsp}, namely $*$-algebras, denoted later by $\mathcal{O}_n^q$,
which are generated by $a_i$, $a_i^*$, $i=1$,~\dots,~$n$, subject to
the relations
\begin{equation}\label{qqq}
a_i^*a_i=1+q a_ia_i^*,\quad a_i^*a_j=0,\quad i\ne j,\
i,j=1,\ldots,n,\ 0<q<1.
\end{equation}
Note that  for $q=0$ we obtain the  $*$-algebra,
$\mathcal{O}_n^{0}$,  generated by $n$ isometries with orthogonal
ranges,  i.e.
\[
\mathcal{O}_n^{0} =\mathbb{C}\left< s_i,\ s_i^*s_j=\delta_{ij}1,\
  i,j=1,\ldots,n \right>.
\]
Its enveloping $C^*$-algebra is an extension of the Cuntz
$C^*$-algebra $O_n$ by compact operators whose representation
theory  was extensively studied (see \cite{bra_jor_ost04, cun77}
and references therein).

 Our aim is to describe irreducible unbounded representations
of $\mathcal{O}_n^q$, $n>1$. Note that each $O_n^q$, $n>1$, has also
bounded representations, however, since the corresponding universal
enveloping $C^*$-algebra $O_n^q$  is isomorphic to the
Cuntz-Toeplitz $C^*$-algebra ${O}_n^0$ (see Section 1),
$\mathcal{O}_n^q$ is not of type I algebra and the problem of
unitary classification of its irreducible $*$-representations is
complicated. Nevertheless it turns out to be possible to classify
all its ``well-behaved'' unbounded irreducible $*$-representations
up to  unitary equivalence. Unbounded representations are known to
be a very delicate thing, since  unbounded operators are not defined
on the whole space. Depending on  chosen domains of representations
they
 can behave differently
(see, e.g. \cite{nel59,sch90}). In the theory of $*$-representations of
finite-dimensional Lie algebras the class of well-behaved
representations form the representations which can be integrated
to unitary representations of the corresponding simply connected
Lie group. Nelson's fundamental theorem  (see
\cite{barut_raczka,nel59}) gives a criterion for the integrability
in terms of the Laplace operator of the Lie algebra, requiring its
essential self-adjointness on a common invariant dense domain. Our
definition of {\it well-behaved} representation is motivated by
this issue.

The  well-known Stone-von Neumann theorem says that up to unitary
equivalence there exists a unique irreducible representation of CCR
which is unbounded and given by raising and lowering operators; the
representation  is often called the  Fock representation. However,
for $q$-CCR and, as we will see, for  ${\mathcal O}_n^q$, $q\in
(0,1)$, the Fock representation is bounded, and a whole bunch of
irreducible unbounded representations  which do not have any
classical analogs, arises (for $q$-CCR see for example
\cite{ost_sam89romp}).

The paper is organized as follows. In Preliminaries we recall a
classification of irreducible
representations of $\mathcal{O}_1^q$ (or $q$-CCR) and  prove an
isomorphism of the enveloping $C^*$-algebras of $\mathcal{O}_n^q$,
$q\in (0,1)$, to that of ${\mathcal O}_n^0$. As a consequence we get that
the description of
bounded $*$-representations of $\mathcal{O}_n^q$ is equivalent to
that  of $*$-representations of the Cuntz-Toeplitz algebra ${\mathcal O}_n^0$.

In Section 2.1 we give a definition of well-behaved unbounded
$*$-representation of $\mathcal{O}_n^q$ in spirit of
\cite{pus_wor89} and also present an equivalent one in terms of
bounded operators.

Finally, in Section 2.2 we obtain a classification of all irreducible
well-behaved unbounded
$*$-representations of $\mathcal{O}_n^q$ up to unitary
equivalence.

\section{Preliminaries}
Recall some facts on representation theory and properties of
the universal enveloping $C^*$-algebra of $q$-CCR or $\mathcal{O}_1^q$, $0\le
q<1$, see \cite{ost_sam99}. If $q=0$ we get the well-known
$*$-algebra generated by a single isometry. Obviously any
representation of $\mathcal{O}_1^0$ is bounded. Moreover, any
irreducible representation of ${\mathcal O}_{1}^0$, up to a
unitary equivalence, is either one-dimensional
$a_1=\exp{\imath\varphi}$, $\varphi\in [0,2\pi)$, or the
infinite-dimensional, called also the {\em Fock} representation,
given by the action
\[
a_1\, e_n=e_{n+1},\quad n\in\mathbb{N},
\]
on an orthonormal basis $\{e_n\colon n\in{\mathbb N}\}$ in
$l_2({\mathbb N})$.

When $0<q<1$, unbounded representations will also arise.
Defining ``well-behaved''  unbounded representations as in
Section~\ref{sec_unbounded} one has the following
\begin{proposition}\label{onedim}
Any irreducible representation of $\mathcal{O}_1^q$ with $0<q<1$
is unitarily equivalent to exactly one listed below.
\\ \noindent
1. {\em The  Fock representation} acting on $l_2(\mathbb{N})$:
\[
\pi_F(a_1)\ e_n=\sqrt{\frac{1-q^n}{1-q}}\ e_{n+1},\quad n\in\mathbb{N}.
\]
\noindent 2. {\em One-dimensional representations:}
\[
\pi_{\varphi}(a_1)\ e = \sqrt{\frac{1}{1-q}}\exp(i\varphi)\
e,\quad \varphi\in [0,2\pi).
\]
\noindent 3. {\em Unbounded representations} acting on
$l_2(\mathbb{Z}):$
\[
\pi_x(a_1)\ e_n=\sqrt{\frac{1-q^n}{1-q}+q^n x}\ e_{n+1},\quad
n\in\mathbb{Z},\ x\in (1+qx_0,x_0]
\]
where $x_0>\frac{1}{1-q}$ is fixed.
\end{proposition}
One of the fundamental facts on  representation theory of
$\mathcal{O}_1$ is the Wold decomposition theorem, stating that any
isometric operator is an orthogonal direct sum of a multiple of the
unilateral shift and a unitary operator (see \cite{murphy}). Using
the description of irreducible representations of $\mathcal{O}_1^q$
one can get a generalization of the Wold decomposition theorem to
the case of linear operator satisfying $q$-canonical commutation
relation (below we will refer to this fact as the $q$-Wold
decomposition theorem).
\begin{theorem}\label{qccr}
Let $A\colon\mathcal{H}\rightarrow\mathcal{H}$ be a bounded linear
operator satisfying for some $q\in (0,1)$ the $q$-commutation
relation
\[
A^* A=1+q A A^*.
\]
Then $\mathcal{H}$ can be decomposed into orthogonal sum of
subspaces $\mathcal{H}=\mathcal{H}_0\oplus\mathcal{H}_u$ invariant
with respect to the actions of $A$, $A^*$ and such that the restriction
$A|_{\mathcal{H}_0}$ is a multiple of the  weighted shift on
$l_2(\mathbb{Z}_{+})$ defined on the standard basis by $A
e_n =\sqrt{\frac{1-q^n}{1-q}}e_{n+1}$, and
$A|_{\mathcal{H}_u}=\frac{1}{\sqrt{1-q}}U$, where $U$ is a unitary
operator on $\mathcal{H}_u$.
\end{theorem}
If we do not assume an operator $A$ satisfying (\ref{q1}) to be
bounded we  can still decompose $\mathcal{H}$ into orthogonal
direct sum of invariant subspaces $\mathcal{H}_0$ and
$\mathcal{H}_u$ such that the restriction of $A$ to
$\mathcal{H}_0$  is a multiple of the Fock representation.
However in general situation $\mathcal{H}_u$ is decomposed into
orthogonal sum $\mathcal{H}_u =\mathcal{H}_1\oplus\mathcal{H}_2$,
where $A_{|\mathcal{H}_1}=\frac{1}{\sqrt{1-q}}U$ with unitary $U$,
and the restriction $A_{|\mathcal{H}_2}$ is unbounded and given by
a direct integral of unbounded irreducible representations. In
particular, in the polar decomposition $A_{|\mathcal{H}_2}=SC$ the
isometric part $S$ is unitary.

Note that the Fock representation $\pi_F$ of $\mathcal{O}_1^q$ is
faithful and  the universal enveloping $C^*$-algebra of
$\mathcal{O}_1^q$  is isomorphic to the $C^*$-algebra generated
by $\pi_F(a_1)$.
\begin{definition}
Let $A$ be a $*$-algebra. Assume that the set, $\rep\mathcal{A}$, of
all its bounded representations is not empty and
\[
\sup_{\pi\in \rep \mathcal{A}} \Vert \pi(a)\Vert <\infty
\]
for any $a\in\mathcal{A}$. The universal enveloping $C^*$-algebra
of a $*$-algebra $\mathcal{A}$ is the completion of
$\mathcal{A}/\mathcal{R}$ with respect to the following norm
\[
\Vert a+\mathcal{R}\Vert=\sup_{\pi\in\rep
\mathcal{A}}\Vert\pi(a)\Vert,
\]
where
\[
\mathcal{R}=\{a\in\mathcal{A}\mid \pi(a)=0,\ \pi\in \rep
\mathcal{A}\}.
\]
\end{definition}
The existence of the  enveloping $C^*$-algebra
 of $\mathcal{O}_n^q$ follows from the fact that $\Vert
\pi(a_i)\Vert^2\le\frac{1}{1-q}$ , $i=1,\ldots,n$ for any bounded representation
$\pi$ of $\mathcal{O}_n^q$ (see Theorem~\ref{qccr}).
In what follows we denote this $C^*$-algebra  by $O_n^q$. It is known that $O_1^q\simeq O_1^0$ for any $q\in
(-1,1)$, see for example
\cite{jor_sch_wer94}. The same is true for $O_n^q$, $n\in\mathbb{N}$.

\begin{theorem}
$O_n^q\simeq O_n^0$ for any $q\in (-1,1)$.
\end{theorem}

\begin{proof}
Suppose that ${O}_n^q$ is realized by linear operators on a
Hilbert space. Consider the polar decompositions $a_i=s_ic_i$,
$i=1,\ldots,d$, where $c_i^2=a_i^*a_i$ and $s_i$ are partial
isometries such that $\ker s_i=\ker c_i$. Since the spectrum
$\sigma(c_i^2)$ of $c_i^2$ is $\{\frac{1-q^n}{1-q},\
n\in\mathbb{N}\}\cup\{\frac{1}{1-q}\}$, each $c_i$ is invertible
and therefore each $s_i$ is an  isometry. Moreover, $s_i=a_i
(a_i^*a_i)^{-\frac{1}{2}}\in C^*(a_i,a_i^*)$, $i=1,\ldots,n$.
Further from $a_i^*a_j=0$ one has
\[
c_is_i^*s_jc_j=0,\quad \mbox{hence}\ s_i^*s_j=0,\ i\ne j.
\]
Since in any irreducible  bounded representation of
$\mathcal{O}_1^q$ with $-1<q<1$ one has
\[
a_1 = s_1 \Bigl(\sum_{n=0}^{\infty}q^n s_1^n s_1^{*n}\Bigr),
\]
the same equality holds in $O_n^q$ for $a_i$ and $s_i$,
$i=1,\ldots,n$.

Therefore $a_i\in C^* (s_i,s_i^*,\; i=1,\ldots,n)$, $i=1,\ldots,n$, giving
 the statement of
the  theorem.\hfill$\Box$

Since the Cuntz-Toeplitz algebra $O_n^0$ is  a not of type I
algebra (see \cite{bra_jor_ost04}), this theorem shows that the
classification problem of all irreducible representations of
$O_n^q$ and therefore all irreducible bounded representations of
${\mathcal O}_n^q$ is very complicated.
\end{proof}

\section{Unbounded representations of
$\mathcal{O}_n^q$ }
\subsection{Definition and properties}\label{sec_unbounded}
We will start with defining well-behaved representations of
$\mathcal{O}_n^q$.
\begin{definition}
We say that a family $\{A_i, i=1,\ldots,n\}$ of closed operators
on a Hilbert space ${\mathcal H}$ defines a well-behaved
representation $\pi$ of $\mathcal{O}_n^q$ if
\begin{enumerate}
\item there exists a dense linear subset ${\mathcal
D}\subset{\mathcal
    H}$ which is invariant with respect to $A_1,\ldots,A_n$,
  $A_1^*,\ldots,A_n^*$;
\item $A_i^*A_ix=(1+qA_iA_i^*)x$, $A_i^*A_jx=0$, $i\ne j$,
  $i=1,\ldots,n$, $x\in {\mathcal D}$;
\item the positive linear operator $\Delta=\sum_{i=1}^n A_i^*A_i$
is
  essentially selfadjoint on ${\mathcal D}$.
\end{enumerate}
\end{definition}
This definition is similar to the definition of unbounded representations of twisted canonical commutation relations given by Pusz and Woronowicz
(\cite{pus_wor89}). It is motivated by the Nelson criterion of the integrability for representations of Lie algebras (\cite{barut_raczka,nel59}).

Next two theorems provide a criteria for representations to be well-behaved
 in terms of bounded operators. It will be important for later classification of well-behaved (irreducible) representations of ${\mathcal O}_n^q$.
Before stating the theorems we  recall the definition of analytic and bounded
vectors for an
unbounded operators, which will be used in the proofs.

If $A$ is an operator in a Hilbert space $H$, then $u\in H$ is said to be an {\it analytic vector}  (a {\it bounded vector}) for $A$, if
$$\sum_{n=0}^{\infty}\frac{||A^nu||}{n!}s^n<\infty,$$
for some $s>0$ ($||A^nu||\leq C^n$, for some $C>0$ respectively).
These concepts are fundamental in the theory of integrable representations of Lie algebras (\cite{barut_raczka, fs3,nel59}).

\begin{theorem}\label{domain}
 Let $\{A_i, i=1,\ldots,n\}$ be a family of linear operators on
 ${\mathcal H}$ defining a well-behaved representation of $\mathcal{O}_n^q$ and  let $A_i=S_i|A_i|$ be the
 polar decomposition of  $A_i$ and $D_i=S_i|A_i|S_i^*$, $i=1,\ldots,n$. Then

(a) $f(D_i^2)S_i=S_if(1+qD_i^2)$  for any real bounded measurable
function $f$;

(b) for any $i$, $j$ the operators $|A_i|$, $|A_j|$ commute in the
sense of resolutions of identity;

(c) $S_i^*S_j=\delta_{ij}I$.
\end{theorem}

\begin{proof}
Let $C_j^2$ be the closure of $A_j^*A_j$ on ${\mathcal D}$,
$j=1,\ldots,n$. Clearly, $C_j^2$ is symmetric. In order to show
that all $C_j^2$ are selfadjoint and mutually commute in the sense
of resolution of identity we  prove first that
\begin{equation}\label{f1}
\Delta^nC_j^2y=C_j^2\Delta^ny, y\in {\mathcal D}(\Delta^{n+2}).
\end{equation}
Here and subsequently, ${\mathcal D}(a)$ denotes the domain of an
operator $a$.

 We have
\begin{eqnarray*}\label{f2}
C_i^2C_j^2&=&(1+qA_iA_i^*)(1+qA_jA_j^*)=
(1+qA_iA_i^*+qA_jA_j^*)\\&=&(1+qA_j^*A_j)(1+qA_iA_i^*)
=C_j^2C_i^2\nonumber
\end{eqnarray*}
on ${\mathcal D}$. Thus if $x$, $y\in{\mathcal D}$ then
\begin{equation}\label{com}
(\Delta x, C_j^2y)=(C_j^2x,\Delta y) \mbox{ and } (C_j^2\Delta x,
y)=(C_j^2x,\Delta y).
\end{equation}
As ${\mathcal D}$ is a core for $\Delta$ the second equality holds
also for any $y\in {\mathcal D}(\Delta)$.

We shall show next that ${\mathcal D}(C_j^2)\supset{\mathcal
D}(\Delta^2)$. In fact, by (\ref{f2})
$$\Delta^2=(\sum_{i=1}^n C_i^2)(\sum_{i=1}^n C_i^2)=\sum_{i=1}^n C_i^4+
\sum_{i\ne j}(1+qA_iA_i^*+qA_jA_j^*)$$ giving $\Delta^2\geq C_j^4$
on ${\mathcal D}$. Since ${\mathcal D}$ is a core for $\Delta^2$
we have that for $y\in{\mathcal D}(\Delta^2)$ there exists
$\{y_n\}\in{\mathcal D}$ such that $y_n\to y$, $\Delta^2 y_n\to
\Delta^2 y$ and
\begin{eqnarray*}
||C_j^2(y_n-y_m)||&=&(C_j^2(y_n-y_m),C_j^2(y_n-y_m))
\\
&=&(C_j^4(y_n-y_m),y_n-y_m)
\\
&\leq& (\Delta^2(y_n-y_m),y_n-y_m).
\end{eqnarray*}
Thus the sequence $\{C_j^2y_n\}$ converges to some $z\in{\mathcal
H}$ and by the closededness of $C_j^2$, $y\in {\mathcal D}(C_j^2)$
and $z=C_j^2y$. Moreover,
\begin{equation}\label{delta2}
C_j^4\leq\Delta^2 \mbox{ on }{\mathcal D}(\Delta^2).
\end{equation}

Let now
 $y\in {\mathcal D}(\Delta^3)$. Then $\Delta y\in{\mathcal
   D}(\Delta^2)\subset {\mathcal D}(C_j^2)$, $y\in{\mathcal
   D}(\Delta^2)\subset{\mathcal D}(C_j^2)$ and by
(\ref{com})
$$(\Delta x,C_j^2y)=(x,C_j^2\Delta y), x\in{\mathcal D}.$$
By the closededness argument the same equality holds for any
$x\in{\mathcal D}(\Delta)$ giving $C_j^2{\mathcal
  D}(\Delta^3)\subset{\mathcal D}(\Delta)$ and
$\Delta C_j^2=C_j^2\Delta$ on ${\mathcal D}(\Delta^3)$. We proceed
now by induction and suppose that for any $y\in {\mathcal
D}(\Delta^n)$, $n\geq 3$, $C_j^2y\in{\mathcal D}(\Delta^{n-2})$
and
\begin{equation}\label{ind_step}
\Delta^{n-2}C_j^2y=C_j^2\Delta^{n-2}y.
\end{equation}
In particular, if $y\in{\mathcal D}(\Delta^{n+1})$ then
$C_j^2y\in{\mathcal D}(\Delta^{n-2})$ and (\ref{ind_step}) holds.
Let $z=\Delta^{n-2}y$. Then $z\in {\mathcal D}(\Delta^3)$ and
$C_j^2z=\Delta^{n-2}C_j^2y\in{\mathcal D}(\Delta)$. Therefore
$C_j^2y\in{\mathcal D}(\Delta^{n-1})$ and
$$
\Delta^{n-1}C_j^2y=\Delta\Delta^{n-2}C_j^2y=\Delta
C_j^2\Delta^{n-2}y=C_j^2\Delta^{n-1}y.
$$
Then by induction we obtain (\ref{f1}) for any $n\geq 1$.

Let ${\mathcal D}_{\omega}$ be the space of analytic vectors for
$\Delta$. As ${\mathcal D}_{\omega}\subset {\mathcal
D}(\Delta^k)$, for any $k\in{\mathbb N}$, by (\ref{delta2}) and
(\ref{ind_step})
\begin{eqnarray*}
||\Delta^nC_j^2x||^2&=&||C_j^2\Delta^nx||^2=(C_j^2\Delta^nx,C_j^2\Delta^nx)=
(\Delta^nC_j^4\Delta^nx,x)\\
&\leq& (\Delta^n\Delta^2\Delta^nx,x)=||\Delta^{n+1}x||^2.
\end{eqnarray*}
 This shows that $C_j^2{\mathcal D}_{\omega}\subset{\mathcal
D}_{\omega}$ and, moreover, $C_j^2$ mutually commute on ${\mathcal
D}_{\omega}$. The last can be seen by computing the scalar product
of $C_i^2C_j^2x-C_j^2C_i^2x$, $x\in{\mathcal D}$, with $y\in
{\mathcal D}_{\omega}$.

Next we show that any $x\in {\mathcal D}_{\omega}$ is also
analytic for all $C_i^2$. In fact, as $C_i^4\leq \Delta^2$ on
${\mathcal
  D}_{\omega}$, by assuming by induction that
$C_i^{4n}\leq\Delta^{2n}$ on ${\mathcal D}_{\omega}$ we obtain
\begin{eqnarray*}
C_i^{4(n+1)}&=&C_i^2C_i^{4n}C_i^2\leq C_i^2\Delta^{4n}C_i^2
\\
&=&
\Delta^{2n}C_i^4\Delta^{2n}\leq\Delta^{2n}\Delta^2\Delta^{2n}=\Delta^{2(n+1)}.
\end{eqnarray*}
This gives
$$||(C_i^2)^nx||^2=(x,C_i^{4n}x)\leq (x,\Delta^{2n} x)\leq ||\Delta^nx||^2,$$
i.e. ${\mathcal D}_{\omega}$ is a subset of analytic vectors for
all $C_i^2$. As ${\mathcal D}_{w}$ is invariant with respect to
all $C_i^2$ and $C_i^2$ mutually commute on ${\mathcal D}_{\omega}$, we have that all
$C_i^2$ are selfadjoint and  mutually  strongly commute, i.e. in
the sense of resolutions of  identity. In particular, we have
proved that each $C_i^2$ is essentially selfadjoint on ${\mathcal
D}$, and $C_i^2=A_i^*A_i=|A_i|^2$.

Next we prove that $f(C_i^2)S_i=S_if(1+qC_i^2)$ for any bounded
measurable function $f$. Let $R_i=C_i^2$ for notation simplicity.
As $R_iA_i=A_i(1+qR_i)$ and $R_iA_i^*=A_i^*(R_i-1)/q$ on
${\mathcal D}$, using arguments similar to one given above one can
show that for any non-negative integer $n$ and $x\in {\mathcal
D}(R_i^{k-1})$ we have that $A_i^*x\in{\mathcal D}(R_i^{k})$,
$A_ix\in{\mathcal D}(R_i^{k})$ and
$$R_i^kA_ix=A_i(1+qR_i)^kx \mbox{ and }R_i^kA_i^*x=A_i^*((R_i-1)/q)^kx.$$

Taking now ${\mathcal D}_{i,\omega}$ the space of analytic vectors
for $R_i$ and using that $R_i\leq (1+R_i)^2$ we obtain
\begin{eqnarray*}
||R_i^kA_ix||^2&=&||A_i(1+qR_i)^kx||^2=((1+qR_i)^kA_i^*A_i(1+qR_i)^kx,x)\\
&\leq&
((1+qR_i)^k(1+R_i)^2(1+qR_i)^kx,x)\\ &=&((1+R_i)(1+qR_i)^{2k}(1+R_i)x,x)\\
&\leq& ((1+R_i)^{2(k+1)}x,x)=||(1+R_i)^kx||^2
\end{eqnarray*}
giving that $A_i{\mathcal D}_{i,\omega}\subset {\mathcal
D}_{i,\omega}$. Analogously, one proves that $A_i^*{\mathcal
D}_{i,\omega}\subset {\mathcal D}_{i,\omega}$. Moreover, the
relations $A_i^*A_i=1+qA_iA_i^*$  and $R_iA_i=A_i(1+qR_i)$ hold on
${\mathcal D}_{i,\omega}$ which can be shown analogously to
commutation of $C_i^2$ on ${\mathcal D}_{\omega}$ above. That
${\mathcal D}_{i,\omega}$ is a core for $A_i$ and $A_i^*$ can be
proved using the arguments in \cite[Proposition 3.3]{pus_wor89}.
The condition $f(C_i^2)S_i=S_if(1+qC_i^2)$   now follows from
\cite[Theorem~1, Theorem~2]{ost_sam89romp}. Furthermore,
\begin{eqnarray}\label{isom}
A_i^*A_i=1+qA_iA_i^*\Leftrightarrow C_i^2=1+qS_iC_i^2S_i^*\Rightarrow\nonumber\\
S_iS_i^*C_i^2=S_iS_i^*+qS_iC_i^2S_i^*=C_i^2-1
\end{eqnarray}
 on ${\mathcal D}$
giving
$$(1-S_iS_i^*)C_i^2=(1-S_iS_i^*)$$
As ${\mathcal D}$ is a core for $C_i^2$ we obtain the equality on
${\mathcal D}(C_i^2)$. Similarly,
$$(1-S_iS_i^*)C_i^{2n}=(1-S_iS_i^*)$$ on ${\mathcal D}(C_i^{2n})$ and
in particular on ${\mathcal D}_{i,\omega}$ for any $n\geq 1$. The
arguments similar to one in \cite[Theorem~1]{ost_sam89romp} give
$(1-S_iS_i^*)f(C_i^2)=f(1)(1-S_iS_i^*)$ for any bounded Borel
function $f$. This will also give that  $S_iS_i^*$ commute with
resolution of the identity of $C_i$ and from (\ref{isom}) we will
get that $C_i\geq 1$ and since $\ker C_i=\ker S_i$,
 $S_i$ is an isometry.

To obtain $f(D_i^2)S_i=S_if(1+qD_i^2)$ we note that
$A_iA_i^*=D_i^2$ and $(A_iA_i^*)A_i=A_i(1+qA_iA_i^*)$ on
${\mathcal D}_{i\omega}$. Moreover, clearly, any vector in
${\mathcal D}_{i,\omega}$ is analytic for $A_iA_i^*$. Using again
 \cite[Theorem~1, Theorem~2]{ost_sam89romp} we obtain the desired equality.

What is left to prove is that $S_i^*S_j=0$ if $i\ne j$. We have
$A_i^*A_j=C_iS_i^*S_jC_j=0$ on ${\mathcal D}$ and
$(S_i^*S_jC_jx,C_iy)=0$ for any $x\in{\mathcal D}$, $y\in{\mathcal
D}(C_i)$. As each $S_i$ is an isometry, the range of $C_i$ is
dense implying $S_i^*S_j=0$.\hfill$\Box$
\end{proof}
\begin{remark}\rm
Let $E_j(\cdot)$ be the resolution of identity of
$D_j^2=S_j|A_j|^2S_j^*$, $j=1,\ldots,n$. Then, by \cite{ost_sam89romp}, $(a)$ is equivalent to
\begin{equation}\label{f3}
E_j(\delta)S_j=S_jE_j(q^{-1}(\delta-1)), j=1,\ldots,n \mbox{ for
any Borel }\delta\subset\mathbb{R}. \end{equation}
\end{remark}

We will use the following notation when proving our next theorem.

Let $\Lambda=\{\emptyset,\
\alpha=(\alpha_1,\alpha_2,\ldots,\alpha_k),\ 1\le\alpha_i\le n\,\
k\in\mathbb{N}\}$ be the set of all finite multi-indices. Introduce
the transformations $\sigma,\
\sigma_k\colon\Lambda\rightarrow\Lambda$, $k=1,\ldots,n$
\begin{eqnarray*}
\sigma(\alpha_1,\ldots,\alpha_s)& =&(\alpha_2,\ldots,\alpha_s),\
\sigma(\alpha_1)=\emptyset,\\
\sigma_k(\alpha_1,\ldots,\alpha_s)&=&(k,\alpha_1,\ldots,\alpha_s).
\end{eqnarray*}

Below having any family $u_1,\ldots,u_n$ of elements of some algebra
and any nonempty multi-index
$\alpha=(\alpha_1,\ldots,\alpha_k)\in\Lambda$ we will denote by
$u_{\alpha}$ the product $u_{\alpha_1}u_{\alpha_2}\cdots
u_{\alpha_k}$, it will be also convenient for us to put
$u_{\emptyset}=1$

\begin{theorem}\label{bounded}
Let $A_i=S_i|A|_i$, $i=1,\ldots, n$, be the polar decompositions of closed operators $A_i$.
If $|A_i|$, $D_i$, $S_i$, $i=1,\ldots, n$ satisfy conditions
(a)-(c) of { Theorem \ref{domain}}, then $\{A_i,
i=1,\ldots, n\}$ defines a well-behaved representation of
${\mathcal O}_n^q$.
\end{theorem}
\begin{proof}
We construct the necessary domain.

The condition (\ref{f3})  implies that given a fixed $j$,  the
operators $A_j$, $A_j^*$ form a (well-behaved) representation of
$q$-CCR relation with one degree of freedom. From the generalized
Wold decomposition for representations of one-dimensional $q$-CCR
we can write $\mathcal{H}=\mathcal{H}_0^{(j)}\oplus
\mathcal{H}_s^{(j)} \oplus \mathcal{H}_u^{(j)}$ so that in
$\mathcal{H}_0^{(j)}$ we have a multiple of the Fock
representation, in $\mathcal{H}_s^{(j)}$ we have
$D_j^2=(1-q)^{-1}$ and in $\mathcal{H}_u^{(j)}$ is unbounded
component. Notice that in $\mathcal{H}_s^{(j)}$ and in
$\mathcal{H}_u^{(j)}$ the operator $S_j$ is unitary. Let
$\mathcal{H}_j$ be a span of the vectors $S_\alpha x$, $x\in
\mathcal{H}_u^{(j)}$, $\alpha \in \Lambda$.
\begin{lemma}\label{decomp}
1. $\mathcal{H}_j$, $j=1$, \dots, $n$  are invariant subspaces.

2. $\mathcal{H}_j\perp \mathcal{H}_k$, $j\ne k$.

3. In the subspace
$\mathcal{H}_0=\mathcal{H}\ominus(\mathcal{H}_1\oplus\dots\oplus
\mathcal{H}_n)$ the representation is bounded.
\end{lemma}

\begin{proof}
1. It is sufficient to show that each $\mathcal{H}_j$ is invariant
with respect to $S_k$, $S_k^*$, and $E_k(\delta)$, $k=1$, \dots, $n$
for any measurable $\delta$.

From $S_kS_\alpha=S_{\sigma_k(\alpha)}$ obviously follows the
invariance with respect to $S_k$, $k=1$, \dots, $n$. For $\alpha\ne
\emptyset$ the invariance with respect to $S_k^*$ follows as well
since $S_k^*S_\alpha=\delta_{k\alpha_1}S_{\sigma(\alpha)}$. For the
vectors of the form $x\in \mathcal{H}_u^{(j)}$, since $S_j$ is
unitary in $ \mathcal{H}_u^{(j)}$, we have $x=S_jS_j^*x$ and
therefore, $S_k^*x=\delta_{jk}S_j^*x\in \mathcal{H}_u^{(j)}$.

Take a measurable $\delta$. Since
 $S_kS_k^*$ is the projection on the cokernel of $D_k^2$, then for
 $\delta$ not containing $\{0\}$ we have
 $E_k(\delta)S_j=E_k(\delta)S_kS_k^*S_j=0$, and
 $E_k(\{0\})s_j=(1-S_kS_k^*)S_j =S_j$, $j\ne k$. Therefore, for
 $0\in\delta$ we have $E_k(\delta)S_j=S_j$, and for $0\notin \delta$
 we have $E_k(\delta)S_j=0$, $j\ne k$. Thus we have that $E_k(\delta)
 S_\alpha x\in \mathcal{H}_j$ for $\alpha_1\ne k$.

If $\alpha_1=k$, then we apply (\ref{f3}) and get
$E_k(\delta)S_\alpha x = S_k E_k(q^{-1}(\delta -1))
S_{\sigma(\alpha)}x \in \mathcal{H}_j$ by induction.

It remains to consider the case $\alpha=\emptyset$. Again, we can
write $x=S_jS_j^*x$, and  apply the arguments above.

2. Take arbitrary $x\in \mathcal{H}_u^{(j)}$, $y \in
\mathcal{H}_u^{(k)}$, $j\ne k$.

i) Since $x=S_jS_j^*x$, $y=S_kS_k^*y$, we have $(x,y)=(S_j^*x,
S_j^*S_kS_k^*y)=0$.

ii) For any $\alpha=(\alpha_1,\dots,\alpha_l)$ we have $(S_\alpha
x,y) = (x, S_\alpha^* y)=(S_j^*x,S_j^*S_\alpha^*y)$. But since
$y=S_k^{l+1}(S^*)^{l+1} y$ and $S_j^*S_\alpha s_k^{l+1}=0$, the
latter scalar product is zero.

iii) For any $\alpha$, $\beta\in \Lambda$ we have $(S_\alpha
x,S_\beta y) =0$ using quite similar arguments.

3) Is obvious, since unbounded component of any $A_j$ in its Wold
decomposition generates $\mathcal{H}_j$.\hfill$\Box$
\end{proof}

Let us continue the proof of the theorem. By { Lemma \ref{decomp}}
we decompose the representation space,
$\mathcal{H}=\mathcal{H}_0\oplus \mathcal{H}_1\dots \oplus
\mathcal{H}_n$, in each component we construct the necessary
domain $\mathcal D_j$, $j=0$, \dots, $n$, and take $\mathcal
D=\mathcal D_0\oplus \mathcal D_1\oplus \dots\oplus \mathcal D_n$.

Since in $\mathcal{H}_0$ the operators are bounded, we can take
$\mathcal D_0=\mathcal{H}_0$.

We fix some $j=1,\ldots,n$ and construct the corresponding set
$\mathcal D_j \subset \mathcal{H}_j$.

Let  $\mathcal D_j$ be a span of vectors of the form $S_\alpha
E_j(\delta) x$, where    $x \in \mathcal{H}_u^{(j)}$, $\alpha \in
\Lambda$, $\delta\subset \mathbb R_+$  are bounded measurable
sets.

1. $\mathcal D_j$ is dense in $\mathcal{H}_j$. Indeed, choose
$\delta=[0,t]$, then for any $x\in \mathcal{H}_u^{(j)}$, $\mathcal
D_j$ contains $x_t=E_j([0,t]) x$, which converges to $x$ strongly
as $t\to\infty$. Applying the operators $S_\alpha$, $\alpha \in
\Lambda$, we obtain total in $\mathcal{H}_j$ set.

2. $\mathcal D_j \subset \mathcal{D}(D_k^2)$, $k=1$, \dots, $n$,
and consists of bounded vectors for these operators.

First we show that for any $z\in\mathcal D_j$, the sequence
$E_k([0,t])D_k^2 z$ converges in $\mathcal{H}_j$ as $t\to\infty$.
But as noticed above $E_k([0,t])S_l = E_k(\{0\})S_l$, $k\ne l$,
and hence $E_k([0,t])D_k^2 S_l=0$, $k\ne l$. Therefore,
$E_k([0,t])D_k^2z=0$ on any $z=S_\alpha E_j(\delta) x$, where
$\alpha=(\alpha_1,\dots,\alpha_m)$ with $\alpha_1\ne k$.

For $\alpha_1=k$ we have $E_k([0,t])D_k^2z=S_k(1+qD_k^2)
E_k([0,q^{-1}(t-1)])S_{\sigma(\alpha)} E_j(\delta)x$ and for
$\alpha$ with $\alpha_1=\dots=\alpha_m=k$
\begin{eqnarray}
&\displaystyle E_k([0,t])D_k^2 S_k^m S_{\sigma^m(\alpha)}E_j(\delta) x \nonumber
\\
&\displaystyle= S_k^m (\frac{1-q^m}{1-q} I +q^m D^2_k)
E_k([0,\frac{(1-q)t-1+q^m}{(1-q)q^m}]) S_{\sigma^m(\alpha)}
E_j(\delta) x.
 \label{kkk}
\end{eqnarray}
If $\sigma^m(\alpha)\ne\emptyset$ and $\alpha_{m+1}\ne k$ the
expression in (\ref{kkk}) is equal to
\[
\frac{1-q^m}{1-q} S_k^m S_{\sigma^m(\alpha)} E_j(\delta) x =
\frac{1-q^m}{1-q} S_\alpha E_j(\delta) x.
\]
Here we use the equalities $D_k^2E_k(\{0\})=0$ and $E_k(\{0\})
S_j=S_j$. Similarly, in the case $\sigma^m(\alpha)=\emptyset$ for
$k\ne j$ using $x=S_jS_j^*x$ we have
\[
E_k([0,t])D_k^2 s_k^m E_j(\delta) x = \frac{1-q^m}{1-q} S_k^m
E_j(\delta)  x.
\]

For $k=j$, $\sigma^m(\alpha)=\emptyset$, and large $t$ we have
$E_j(\delta)E_j([0,\frac{(1-q)t-1+q^m}{(1-q)q^m}])=E_j(\delta)$
since $\delta_j$ is bounded, therefore in this case (\ref{kkk})
does not depend on $t$ as well as above and obviously converges in
$\mathcal{H}_j$ as $t\to\infty$.

Finally, for a bounded $\delta$ and $x \in \mathcal{H}_u^{(j)}$ we
have the following expressions:
\begin{eqnarray}
D_k^2 S_\alpha E_j(\delta) x&=&0, \quad \alpha_1\ne k; \nonumber
\\
D_k^2 S_{\alpha} E_j(\delta) x&=&\frac{1-q^m}{1-q} S_\alpha
E_j(\delta)x, \nonumber
\\
&\qquad&
\alpha=(\underbrace{k,\dots,k}_m,\alpha_{m+1},\dots,\alpha_l),\
\alpha_{m+1}\ne k; \nonumber
\\
D_k^2 S_k^m E_j(\delta) x&=&\frac{1-q^m}{1-q} S_k^mE_j(\delta)x,
\quad k\ne j; \nonumber
\\
D_j^2 S_j^m E_j(\delta) x&=&S_j^mE_j(\delta)  (\frac{1-q^m}{1-q} I
+q^m D^2_j)E_j(\delta) x.\label{d_j},
\end{eqnarray}
where in the last formula we used
$$
E_j(\delta)  (\frac{1-q^m}{1-q}
I +q^m D^2_j)E_j(\delta)=  (\frac{1-q^m}{1-q} I +q^m
D^2_j)E_j(\delta).
$$
 From (\ref{d_j}) we conclude that $D_k$, $k\ne
j$, is bounded in $\mathcal{H}_j$ with $||D_k^2|_{\mathcal{H}_j}||
= (1-q)^{-1}$. For $k=j$ we have $||D_j^{2p}S_\alpha E_j(\delta)
x|| \le (\frac{1}{1-q} + r)^p \|E_j(\delta)x\|$, which means that
$\mathcal D_j$ consists of bounded vectors of $D_j^2$.

3. $\mathcal D_j$ is invariant with respect to $S_k$, $S_k^*$,
$k=1$, \dots, $n$. The invariance with respect to $S_k$ is obvious.
For $z=S_\alpha E_j(\delta) x$ with $x\in \mathcal{H}_u^{(j)}$ we
have $S_k^* z =\delta _{k\alpha_1} S_{\sigma(\alpha)} E_j(\delta) x$
if $\alpha\ne \emptyset$. For $z= E_j(\delta)  x$ we have
$z=S_jS_j^* z$ and
\[
S_k z=\delta_{kj} S_j^* z = \delta_{kj} E_j(q^{-1} (\delta
-1))S_j^*x \in \mathcal D_j.
\]

4. Define $A_k$, $k=1$, \dots, $n$ as a closure of $D_kS_k$ from
$\mathcal D$. Then $\mathcal D \subset D(A_k), D(A_k^*)$, $k=1$,
\dots, $n$ and (\ref{qqq}) holds on $\mathcal D$. This follows
directly.

5. One can easily see that $\mathcal D$ consists of bounded
vectors for the operators $S_j^*D_j^2S_j$ and using their
commutation, for $\Delta$ as well. Therefore $\Delta$ is
essentially selfadjoint on ${\mathcal D}$.\hfill$\Box$
\end{proof}

\begin{remark}
In fact it follows from unitarity of $S_j$ on
$\mathcal{H}_u^{(j)}$ that  $\mathcal{H}_j$ coincides with the
closure of the span of the family $\{S_{\alpha}x,\
x\in\mathcal{H}_u^{(j)},\ \alpha\in\Lambda_j \}$, where
$\Lambda_j=\{\emptyset,\ (\alpha_1,\ldots,\alpha_k),\
1\le\alpha_i\le n,\ \alpha_k\ne j,\ k\in\mathbb{N}\}$.
\end{remark}

\begin{remark}\label{remark3}
It follows from the considerations above that
$S_{\alpha}\mathcal{H}_u^{(j)}$ are invariant with respect to
$D_k^2$ and if $\alpha\ne\emptyset$, $\alpha\in\Lambda_j$, then
restriction of $D_j^2$ to $S_{\alpha}\mathcal{H}_u^{(j)}$ is
bounded. In fact, for nonempty $\alpha\in\Lambda_j$ and any
$k=1,\ldots,n$ one has
\begin{equation}\label{fore}
D_k^2 x=\frac{1-q^{m_k(\alpha)}}{1-q}x,\quad x\in
S_{\alpha}\mathcal{H}_u^{(j)},
\end{equation}
where the function $m_k(\alpha)\in\mathbb{Z}_{+}$,
$\alpha\ne\emptyset$, is determined by the condition
\begin{equation}\label{mk}
\sigma_k^{m_k(\alpha)}\sigma^{m_k(\alpha)}(\alpha)=\alpha,\quad
\sigma_k^{m_k(\alpha)+1}\sigma^{m_k(\alpha)+1}(\alpha)\ne\alpha.
\end{equation}
Recall also that $D_k^2 x=0$, for $x\in\mathcal{H}_u^{(j)}$, $k\ne
j$, and with $m_k(\emptyset)=0$ the formula (\ref{fore}) becomes
true for $x\in\mathcal{H}_u^{(j)}$ and $k\ne j$ also.
\end{remark}
\subsection{Irreducible unbounded representations}
In this section we will give a classification of irreducible
unbounded representations of ${\mathcal O}_n^q$. We will keep
notation from the previous section and consider only well-behaved
representations of ${\mathcal O}_n^q$.

Let $\{ A_1,\ldots, A_n\}$ be a family of closed operators on
${\mathcal H}$ defining a representation $\pi$ of ${\mathcal
O}_n^q$. Let $A_i=S_iC_i$ be the polar decomposition of $A_i$ and
$D_i=S_iC_iS_i^*$. Denote by $E_j(\cdot)$ the resolution of
identity of $D_j$ and let ${\mathfrak B}({\mathbb R})$ be the $\sigma$-algebra of
all Borel subsets of ${\mathbb R}$.
\begin{definition}
The representation $\pi$ will be called {\it irreducible} if the
only operator $C\in B({\mathcal H})$ which commutes with  all
$S_i$, $S_i^*$ and $E_i(\delta)$, $\delta\in {\mathfrak
B}({\mathbb R})$, $i=1,\ldots,n$, is a multiple of unity, or
equivalently, the space ${\mathcal H}$ can not be decomposed into
a direct sum of two non-trivial subspaces which are invariant with
respect to $S_i$, $S_i^*$, and $E_i(\delta)$, $\delta\in
{\mathfrak B}({\mathbb R})$, $i=1,\ldots,n$.
\end{definition}
 It follows from {Lemma \ref{decomp}} that for irreducible unbounded $\pi$ the representation
space ${\mathcal H}$ coincides with
 the closed span ${\mathcal H}_j$ of
$\{S_{\alpha}\mathcal{H}_u^{(j)},\ \alpha\in\Lambda_j\}$ for some
$j\in\{1,\ldots,n\}$.
\begin{lemma}\label{restrict}
Let $\pi$ be an irreducible unbounded representation of
$\mathcal{O}_n^q$ on $\mathcal{H}_j$ for some fixed $j$. Then the
restriction  $A_j|_{\mathcal H_u^{(j)}}$ determines an irreducible
representation of $\mathcal{O}_1^q$ on $\mathcal H_u^{(j)}$.
\end{lemma}

\begin{proof}
Assume that $\mathcal{H}_u^j={H}_u^1\oplus {H}_u^2$ where
${H}_u^i$, $i=1,2$, are invariant with respect to $S_j$, $S_j^*$,
and $E_j(\delta)$ for any Borel $\delta\in\mathbb{R}$. Let
$\mathcal{H}_i$, $i=1,2$, be the closed linear span of
$\{S_{\alpha}{H}_u^i,\ \alpha\in\Lambda_j\}$. Then following
arguments in the proof of { Theorem \ref{bounded}} we obtain
that $\mathcal{H}={H}_1\oplus {H}_2$ where each of ${H}_i$ is
invariant with respect to $S_k$, $S_k^*$, $E_k(\delta)$ for any
$k=1,\ldots,n$ and  $\delta\in\mathbb{R}\in\mathfrak B({\mathbb R})$ contradicting the
irreducibility of $\pi$.\hfill$\Box$
\end{proof}

\begin{theorem}\label{descrep}
Let $j\in\{1,\ldots,n\}$ and $x\in (1+qx_0,x_0]$, where
$x_0>\frac{1}{1-q}$ is fixed. Let $\mathcal{H}$ be a Hilbert space
with orthonormal basis $\{e_{\alpha}^s,\ s\in\mathbb{Z},\
\alpha\in\Lambda_j\}$ and let $A_k$, $k=1,\ldots,n$  be linear
operators  on $\mathcal{H}$ given by
\begin{eqnarray}\label{description}
A_k e_{\alpha}^s & =&
\sqrt{\frac{1-q^{m_k(\sigma_k(\alpha))}}{1-q}}\
e_{\sigma_k(\alpha)}^s,\quad k\ne j\nonumber
\\
A_j e_{\alpha}^s & =& \sqrt{\frac{1-q^{m_j(\sigma_j(\alpha))}}{1-q}}\
e_{\sigma_j(\alpha)}^s,\ \alpha\ne\emptyset\nonumber
\\
A_j e_{\emptyset}^s & =& \sqrt{\frac{1-q^{s}}{1-q}+q^s x}\
e_{\emptyset}^{s+1}
\end{eqnarray}
where $m_k(\alpha)$ are defined by (\ref{mk}). Then
$\{A_1,\ldots, A_k\}$ defines an irreducible representation
$\pi_{(x,j)}$ of $\mathcal O_n^q$. Moreover any unbounded
irreducible representation is unitarily equivalent to exactly one
representation $\pi_{(x,j)}$.
\end{theorem}

\begin{proof}
Let $\mathcal{U}$ be the closure of span of $\{e_{\emptyset}^s,\
s\in\mathbb{Z}\}$. Then $\mathcal{U}$ is invariant with respect to
$A_j$, $A_j^*$ and the restrictions of these operators to
$\mathcal{U}$ determine a well-behaved irreducible representation
of $\mathcal{O}_1^q$. Consider the polar decompositions
$A_k=S_kC_k$, and put $D_k=S_k C_k S_k^*$, $k=1,\ldots,n$. Then
\begin{eqnarray*}
D_k^2 e_{\alpha}^s & =&A_kA_k^*e_{\alpha}^s
=\frac{1-q^{m_k(\alpha)}}{1-q}\ e_{\alpha}^s,\
\alpha\in\Lambda_j,\ \alpha\ne\emptyset
\\
D_k^2 e_{\emptyset}^s & =&A_kA_k^* e_{\emptyset}^s =0,\ k\ne j,\
s\in\mathbb{Z}
\end{eqnarray*}
and $S_k e_{\alpha}^s=e_{\sigma_j(\alpha)}^s$ if either $k\ne j$
or $\alpha\ne\emptyset$, and $S_j
e_{\emptyset}^s=e_{\emptyset}^{s+1}$; $S_k^*
e_{\alpha}^s=\delta_{ki_1}e_{\sigma(\alpha)}^s$ if
$\alpha\ne\emptyset$, $S_j^* e_{\emptyset}^s=e_{\emptyset}^{s-1}$
and $S_k^* e_{\emptyset}^s=0$, $k\ne j$. In particular,
$e_{\alpha}^s=S_{\alpha}e_{\emptyset}^s$, for any non-empty
$\alpha\in\Lambda_j$ and $s\in\mathbb{Z}$. It is a routine to verify
that  the conditions of { Theorem \ref{bounded}} are satisfied
and formulas (\ref{description}) determine a well-behaved
representation $\pi_{(x,j)}$ of $\mathcal{O}_n^q$ with ${\mathcal
U}=\mathcal{H}_{u}^{(j)}$ and
$\pi_{(x,j)}(a_i)=A_i$.

Each representation $\pi_{(x,j)}$ is irreducible. In fact, let $C\in
B(\mathcal H)$ be a selfadjoint operator commuting with $S_k$,
$E_k(\delta)$, $k=1,\ldots,n$,  $\delta\in\mathfrak B(\mathbb{R})$.
One can easily see that $\mathcal{H}_{u}^{(j)}=E_j(\Delta_x)$, where
$\Delta_x=\{\frac{1-q^n}{1-q}+q^n x,\ n\in\mathbb{Z}\}$, giving that
$\mathcal{H}_{u}^{(j)}$ is invariant with respect to $C$. Put
$C_{\emptyset}=C_{|\mathcal{H}_{u}^{(j)}}$; since the representation
of $\mathcal{O}_1^q$ defined by the restriction of
$\pi_{(x,j)}(a_j)$ to $\mathcal{H}_{u}^{(j)}$ is irreducible, one
has $C_{\emptyset}=\lambda_{\emptyset}\mathbf{1}$,
$\lambda_{\emptyset}\in\mathbb{C}$. Further, using the commutation
of $C$ with all $S_k$ and induction on length of
$\alpha\in\Lambda_j$, we get that $S_{\alpha}\mathcal{H}_u^{(j)}$ is
invariant with respect to $C$  for any $\alpha\in\Lambda_j$.
Denoting by $C_{\alpha}$ the corresponding restriction we obtain
again by induction $C_{\alpha}=\lambda_{\emptyset}\mathbf{1}$ for
any $\alpha\in\Lambda_j$. As
$\mathcal{H}=\oplus_{\alpha\in\Lambda_j}s_{\alpha}\mathcal{H}_{u}^{(j)}$
we conclude that $C=\lambda_{\emptyset}\mathbf{1}_{\mathcal{H}}$ and
$\pi_{(x,j)}$ is irreducible.

Next we show that $\pi_{(x,j)}$ are non-equivalent
representations. Since $\pi_{(x,j)}(a_k)$, are bounded if $j\ne k$,
we have that $\pi_{(x,j)}$ is not unitarily equivalent to
$\pi_{(y,k)}$ when $j\ne k$.

Let $x,x'\in (1+qx_0,x_0]$, $j\in\{1,\ldots, n\}$.  Suppose that
$\pi_{(x,j)}\simeq\pi_{(x',j)}$. Denote the corresponding
representation spaces by $\mathcal{H}$ and $\mathcal{F}$  and let
$V\colon\mathcal{H}\rightarrow\mathcal{F}$ be a unitary operator
giving the equivalence of the representations. Then $V$ gives the
equivalence of representations of $\mathcal{O}_1^q$ defined by the
actions of $\pi_{(x,j)}(a_j)$ and $\pi_{(x',j)}(a_j)$ on
$\mathcal{H}$ and $\mathcal{F}$ respectively. Consider the
decompositions
$\mathcal{H}=\mathcal{H}_u^{(j)}\oplus(\mathcal{H}_{u}^{(j)})^{\perp}$
and
$\mathcal{F}=\mathcal{F}_u^{(j)}\oplus(\mathcal{F}_{u}^{(j)})^{\perp}$.
Recall that summands in these decompositions are invariant with
respect to  $\pi_{(x,j)}(a_j)$, $\pi_{(x,j)}(a_j^*)$
($\pi_{(x',j)}(a_j)$, $\pi_{(x',j)}(a_j^*)$ respectively), the
restriction, $\pi_1(a_j)$, of $\pi_{(x,j)}(a_j)$ (and $\pi_1' (a_j)$
of $\pi_{(x',j)}(a_j)$) to $(\mathcal{H}_{u}^{(j)})^{\perp}$
($(\mathcal{F}_{u}^{(j)})^{\perp}$ respectively) is bounded and
\[
\pi_2(a_j)=\pi_{(x,j)}(a_j)_{|\mathcal{H}_u^{(j)}}=\pi_x(a),\quad
\pi_2'(a_j)=\pi_{(x',j)}(a_j)_{|\mathcal{F}_u^{(j)}}=\pi_{x'}(a)
\]
Since $\pi_1$ and $\pi_2$ are disjoint and the same are $\pi_1'$
and $\pi_2'$ one has that $V\pi_{(x,j)}(a_j)V^*=\pi_{(x',j)}(a_j)$
implies $V=V_2\oplus V_1$ and $V_2\pi_2(a_j)V_2^*=\pi_2'(a_j)$,
hence by { Proposition \ref{onedim}} we get $x=x'$.

Finally, by { Lemma \ref{restrict}} any irreducible
well-behaved representation $\pi$ of $\mathcal{O}_n^q$ acting on
Hilbert space $\mathcal{H}$ corresponds to some fixed
$j=1,\ldots,n$ such that the restriction of $\pi(a_j)$ to
$\mathcal{H}_{u}^{(j)}$ determines an irreducible well-behaved
representation of $\mathcal{O}_1^q$. Then decomposing
$\mathcal{H}=\mathcal{H}_u^{(j)}\oplus
(\mathcal{H}_u^{(j)})^{\perp}$ and taking a unitary operator $V$
of the form $V=V_1\oplus\mathbf{1}$, where $V_1$ is unitary acting
on $\mathcal{H}_{u}^{(j)}$ such that
$V_1^*\pi(a_j)_{|\mathcal{H}_u^{(j)}}V_1=\pi_x(a)$ for some $x\in
(1+qx_0,x_0)$, we obtain by Remark~\ref{remark3} that  $V$ gives
the unitary equivalence of  $\pi$ to $\pi_{(x,j)}$.\hfill$\Box$
\end{proof}

\begin{remark}\rm
 Using the same arguments we can describe irreducible
representations of $\mathcal{O}_n^q$ such that some of $S_j$ is not
a pure isometry (i.e. its Wold decomposition consists of the unitary
part). In this case we have the following two possibilities: either
the corresponding representation is unbounded as described in {
Theorem~\ref{descrep}}, or it is unitarily equivalent to one
determined by the following formula
\begin{eqnarray}\label{description1}
A_k e_{\alpha} & =& \sqrt{\frac{1-q^{m_k(\sigma_k(\alpha))}}{1-q}}\
e_{\sigma_k(\alpha)},\quad k\ne j\nonumber
\\
A_j e_{\alpha} & =& \sqrt{\frac{1-q^{m_j(\sigma_j(\alpha))}}{1-q}}\
e_{\sigma_j(\alpha)},\ \alpha\ne\emptyset\nonumber
\\
A_j e_{\emptyset}& =&\exp ({2\pi\imath\
\phi_j})\sqrt{\frac{1}{1-q}}\ e_{\emptyset},\ \phi_j\in
[0,1)
\end{eqnarray}
on ${\mathcal H}$ with orthonormal basis $\{e_{\emptyset},\
e_{\alpha},\ \alpha\in\Lambda_j\}$. Representations corresponding
to
 different $j=1,\ldots,n$ or $\phi_j\in [0,1)$   are
non-equivalent.

In particular, for  $q=0$ we get a classification of all
irreducible representations of $\mathcal{O}_n$ such that one of
the generators is not a pure isometry.
\end{remark}

\section*{Acknowledgements}
The paper was written when the first and the second authors were
visiting Chalmers University of Technology and Gothenburg
University, Sweden. Warm hospitality and excellent working
atmosphere are greatly acknowledged. The research  was partially
supported by a grant from the Swedish Royal Academy of Sciences as a
part of the program of cooperation with the former Soviet Union. The
third author was also supported by the Swedish Research Council.

\end{document}